\documentclass[a4paper,12pt]{article}

\title{Sign Consistency of the Generalized Elastic Net Estimator}
\author{Wencan Zhu, Eric Adjakossa, C\'eline L\'evy-Leduc, Nils Tern\`es}

\usepackage[utf8x]{inputenc}
\usepackage{hyperref}
\usepackage{graphicx,ulem}
\usepackage{amsmath,amssymb,amsthm, amscd,a4wide,enumerate,enumitem}
\usepackage[numbers]{natbib}
\usepackage{units}
\usepackage{authblk}
\usepackage[toc,page]{appendix}
\usepackage{listings}
\usepackage{anyfontsize}
\usepackage{libertine}
\usepackage[libertine]{newtxmath}
\usepackage{booktabs}
\usepackage{xcolor}
\usepackage{calrsfs} 
\DeclareMathAlphabet{\pazocal}{OMS}{zplm}{m}{n}
\DeclareMathAlphabet\mbi{OML}{cmm}{b}{it}
\DeclareSymbolFont{boldsymbols}{OMS}{cmsy}{b}{n}
\DeclareSymbolFontAlphabet{\mathbfcal}{boldsymbols}

\DeclareMathOperator*{\argmin}{Argmin}

\newcommand{\R}{\mathbb{R}}

\usepackage{array}
\newcommand{\by}{\mathbf{y}}
\newcommand{\bX}{\mathbf{X}}
\newcommand{\bP}{\mathbb{P}}

\newcommand{\bu}{\mathbf{u}}

\newcommand{\bbeta}{\boldsymbol{\beta}}
\newcommand{\bepsilon}{\boldsymbol{\epsilon}}
\newcommand{\bSigma}{\boldsymbol{\Sigma}}

\newcommand{\norm}[1]{\left\lVert#1\right\rVert}

\everymath{\displaystyle}
\newtheorem{thm}{Theorem}[section]
\theoremstyle{definition}

\theoremstyle{plain}
\newtheorem{lem}[thm]{Lemma}

\graphicspath{{graphics/}}

\theoremstyle{plain}

\begin{document}
\maketitle

\abstract{

In this paper, we propose a novel variable selection approach in the framework of high-dimensional linear models where the columns of the design matrix are highly correlated. It consists in rewriting the initial high-dimensional linear model to remove the correlation between the columns of the design matrix
  and in applying a generalized Elastic Net criterion since it can be seen as an extension of the generalized Lasso.
  The  properties of  our  approach  called gEN (generalized Elastic Net)
  are  investigated  both  from  a  theoretical  and  a  numerical  point  of
view. More precisely, we provide a new condition called  GIC (Generalized Irrepresentable Condition) which generalizes the EIC (Elastic Net Irrepresentable
Condition) of \cite{Jia2010} under which we prove that our estimator can recover the positions of the null and non null entries of the coefficients
when the sample size tends to infinity.
We also assess the performance of our methodology using synthetic data and compare it with alternative approaches. Our numerical experiments show that
our approach improves the variable selection performance in many cases.




}

\ \\


\section{Introduction}
Variable selection has become an important and actively used task for understanding or predicting \textcolor{black}{an outcome} of interest in many fields such as medicine \citep{Lu, gunter2011variable, gu2013bayesian, zhu2020}, social media \citep{Tufekci, lin2016does, tomeny2017geographic}, \textcolor{black}{or finance} \citep{sermpinis2018modelling, Amendola, uniejewski2019understanding}. 
Through decades, \textcolor{black}{numerous} variable selection methods have been developed \textcolor{black}{such as subset selection \citep{draper1998applied}
  or regularization techniques \citep{bickel2006regularization}.}

Subset selection methods achieve sparsity by selecting the best subset of relevant variables using the Akaike information criterion \citep{akaike1998information} or the Bayesian information criterion \citep{schwarz1978estimating} but are shown to be NP-hard and could be unstable in practice \citep{welch1982algorithmic, breiman1996heuristics}. 

The \textcolor{black}{regularized variable selection techniques} have become popular for their capability to overcome the above difficulties \citep{Lasso, Hastie2005, Zou2005, Wu2009}. 
\textcolor{black}{Among them, the Lasso approach \citep{Lasso}
  is one of the most popular and can be defined as follows. Let $\by$ satisfy the following linear model}
\begin{equation}
\label{lm}
    \by=\bX\bbeta^\star+\bepsilon,
\end{equation}
where $\by=(y_{1}, \ldots, y_{n})^{'}\in \R^n$ is the response variable, \textcolor{black}{$'$ denoting the transposition, $\bX=(\bX_{1}, \ldots, \bX_{p})$ is the design matrix with $n$ rows of observations on $p$ covariates, $\bbeta^\star=(\beta_{1}^\star, \ldots, \beta_{p}^\star)^{'}\in \R^p$ is a sparse vector,
  namely contains a lot
  of null components, and $\bepsilon$ is a Gaussian vector with zero-mean and a covariance matrix equal to $\sigma^{2}\mathbb{I}_{n}$, $\mathbb{I}_{n}$ denoting
  the identity matrix in $\R^n$.} The Lasso \textcolor{black}{approach} estimates $\bbeta^\star$ with a
sparsity enforcing constraint by minimizing the following penalized least-squares criterion:
\begin{equation}\label{lasso}
  L^{Lasso}_{\lambda}(\bbeta)=\norm{\by-\bX\bbeta}_{2}^{2}+\lambda\norm{\bbeta}_{1}, 
\end{equation}
 \textcolor{black}{where $\norm{a}_{1}=\sum_{k=1}^p|a_k|$ denotes the $\ell_1$ norm of the vector $(a_1,\dots,a_p)'$,
   $\norm{b}_{2}^2=\sum_{k=1}^n b_k^2$ denotes the $\ell_2$ norm of the vector $(b_1,\dots,b_n)'$}, and $\lambda$ is
 \textcolor{black}{a positive constant corresponding to} the regularization parameter.
 \textcolor{black}{The Lasso popularity largely comes from the fact that the resulting estimator}
$$
\widehat{\bbeta}^{Lasso}(\lambda)=\argmin_{\bbeta\in\R^p}L^{Lasso}_{\lambda}(\bbeta)
$$ 
is sparse (has only a few nonzero entries), and sparse models are often preferred for their \textcolor{black}{interpretability} \citep{Zhao2006}.
\textcolor{black}{Moreover, $\widehat{\bbeta}^{Lasso}(\lambda)$ can be proved to be sign consistent under some assumptions, namely there exists $\lambda$ such that
  $$\lim_{n\to\infty}\bP\left(sign\left(\widehat{\bbeta}^{Lasso}(\lambda)\right)=sign(\bbeta^\star) \right)=1,$$
  where $\textrm{sign}(x)=1$ if $x>0$, -1 if $x<0$ and 0 if $x=0$.}
\textcolor{black}{Before giving the conditions under which \citep{Zhao2006} prove the sign consistency of $\widehat{\bbeta}^{Lasso}$,
  we first introduce some notations.
  Without loss of generality, we shall assume as in \citep{Zhao2006}
  that the first $q$ components of $\bbeta^\star$ are non null (\textit{i.e.} the components that are associated to the active variables, and 
  denoted as $\bbeta_1^\star$)
  and the last $p-q$ components of $\bbeta^\star$ are null (\textit{i.e.} the components that are associated to the non active variables, and 
  denoted as $\bbeta_2^\star$).
Moreover,
we shall denote by $\bX_1$ (resp. $\bX_2$) the first $q$ (resp. the last $p-q$) columns of $\bX$.
Hence, $C_n= n^{-1}\bX'\bX$, which is the empirical covariance matrix of the covariates,  can be rewritten as follows:}
$$
C_n=
    \begin{bmatrix}
   C_{11}^n &  C_{12}^n \\
    C_{21}^n &  C_{22}^n 
    \end{bmatrix},
    $$    
    \textcolor{black}{with $C_{11}^n=n^{-1}\bX_{1}^{'}\bX_{1}$, $C_{12}^n=n^{-1}\bX_{1}^{'}\bX_{2}$, $C_{21}^n=n^{-1}\bX_{2}^{'}\bX_{1}$,
    $C_{22}^n=n^{-1}\bX_{2}^{'}\bX_{2}$.}
\textcolor{black}{It is proved by Zhao and Yu in \citep{Zhao2006} that $\widehat{\bbeta}^{Lasso}(\lambda)$ is sign
  consistent when the following Irrepresentable Condition (IC) is satisfied:
\begin{equation}\label{IC}
  \left|\left(C_{21}^{n}(C_{11}^{n})^{-1}\textrm{sign}(\bbeta_{1}^\star)\right)_j\right| \leq 1-\alpha, \textrm{ for all } j,
\end{equation}
where $\alpha$ is a positive constant.}
\textcolor{black}{In the case where $p\gg n$, Wainwright develops in \citep{wainwright2009sharp} the necessary and sufficient conditions, for both deterministic and random designs, on $p$, $q$, and $n$ for which it is possible to recover the positions of the null and non null components of $\bbeta^\star$, namely
  its support, using the Lasso.}   

    When there are high correlations between covariates, especially the active ones, the $C_{11}^n$ matrix may not be invertible,
    \textcolor{black}{and the Lasso estimator
    fails to be sign consistent. To circumvent this issue,} Zou and Hastie \citep{Zou2005} introduced the Elastic Net estimator defined \textcolor{black}{by:}
\begin{equation}\label{eq:def_EN}
\widehat{\bbeta}^{EN}(\lambda, \eta)=\argmin_{\bbeta\in\R^p}L^{EN}_{\lambda, \eta}(\bbeta),
\end{equation}
where
$$
L^{EN}_{\lambda, \eta}(\bbeta)=\norm{\by-\bX\bbeta}_{2}^{2}+\lambda\norm{\bbeta}_{1}+\eta\norm{\bbeta}_{2} \textrm{ with }\lambda,\eta>0.
$$
\textcolor{black}{Yuan and Lin prove in \citep{Yuan2007} that when the following Elastic Net Condition (EIC) is satisfied the Elastic Net estimator defined by
  (\ref{eq:def_EN}) is sign consistent when $p$ and $q$ are fixed: there exist positive $\lambda$ and $\eta$ such that} 
\begin{equation}\label{EIC}
  \left|\left(C_{21}^{n}\left(C_{11}^{n}+\dfrac{\eta}{n}\mathbb{I}_{q}\right)^{-1}\left(\textrm{sign}(\bbeta_{1}^\star)+\dfrac{2\eta}{\lambda}\bbeta_1^\star\right)\right)_j\right| \leq 1-\alpha, \textrm{ for all } j.
\end{equation}
\textcolor{black}{Moreover,} when $p$, $q$, and $n$ go to infinity \textcolor{black}{with} $p\gg n$, Jia and Yu \textcolor{black}{prove in \citep{Jia2010} that the sign consistency of the Elastic Net estimator holds if additionally to Condition (\ref{EIC})
$n$ goes to infinity at a rate faster than
$q\log(p-q)$.}


\textcolor{black}{In the case where the active and non active covariates are highly correlated, IC (\ref{IC}) and EIC (\ref{EIC}) may be violated.
  To overcome this issue several approaches were proposed: the Standard PArtial Covariance (SPAC) method \citep{xue2017variable} and preconditioning approaches among others.
Xue and Qu \citep{xue2017variable} developed 
the so-called SPAC-Lasso which enjoys strong sign consistency in both finite-dimensional ($p<n$) and high-dimensional ($p\gg n$) settings. 
However, the authors mentioned that the SPAC-Lasso method only selects the active variables that are not highly correlated to the non active ones, which may
be a weakness of this approach.
The preconditioning approaches consist in transforming the given data 
$\bX$ and $\by$ before applying the Lasso criterion. For example, \citep{jia:2015}
 and \citep{holp:2016} proposed to left-multiply $\bX$, $\by$ and thus
$\bepsilon$ in Model (\ref{lm}) by specific matrices to remove the correlations between
the columns of $\bX$. A major drawback of the latter approach, called HOLP (High
dimensional Ordinary Least squares Projection), is that the preconditioning
step may increase the variance of the error term and thus may alter the
variable selection performance.}

\textcolor{black}{Recently, \cite{zhu2020} proposed another strategy under the following assumption:
  \begin{enumerate}[label=\textbf{\rm (A\arabic*)}]
	\item \label{covX}
  $\bX$ is assumed to be a random design matrix such that its rows  $(\boldsymbol{x}_i)_{1\leq i\leq n}$ are i.i.d. zero-mean Gaussian random vectors
having a covariance matrix equal to  $\boldsymbol{\Sigma}$.
\end{enumerate}
More precisely, they propose to rewrite Model (\ref{lm}) in order to remove the correlation existing between the columns of $\bX$.
Let $\boldsymbol{\Sigma}^{-1/2}:=\boldsymbol{U}\boldsymbol{D}^{-1/2}\boldsymbol{U}^{T}$ where $\boldsymbol{U}$ and $\boldsymbol{D}$ are the matrices involved in the spectral decomposition
of  the symmetric matrix $\boldsymbol{\Sigma}$ given by: $\boldsymbol{\Sigma}=\boldsymbol{U}\boldsymbol{D}\boldsymbol{U}^{T}$, then, denoting
$\widetilde{\bX}=\bX\boldsymbol{\Sigma}^{-1/2}$, (\ref{lm}) can be rewritten as follows:
\begin{equation}\label{lm_modified}
  \by=\widetilde{\bX}\widetilde{\bbeta}^\star+\bepsilon,
\end{equation}
where $\widetilde{\bbeta}^\star=\boldsymbol{\Sigma}^{1/2}\bbeta^\star:=\boldsymbol{U}\boldsymbol{D}^{1/2}\boldsymbol{U}^{T}\bbeta^\star$. 
With such a transformation, the covariance matrix of the $n$ rows of $\widetilde{\bX}$ is equal to identity and the columns of 
$\widetilde{\bX}$ are thus uncorrelated. The advantage of such a transformation with respect to the preconditioning approach proposed by \cite{holp:2016} is that the error term $\bepsilon$ is not modified thus avoiding an increase of the noise which can overwhelm the benefits of a well conditioned design matrix.
Their approach then consists in minimizing the following criterion with respect to $\widetilde{\bbeta}$:
\begin{equation}\label{WLasso}
\norm{\by-\widetilde{\bX}\widetilde{\bbeta}}_{2}^{2}+\lambda\norm{\Sigma^{-1/2}\widetilde{\bbeta}}_{1}, 
\end{equation}
where $\widetilde{\bX}=\bX\boldsymbol{\Sigma}^{-1/2}$ in order to ensure a sparse estimation of $\bbeta^\star$ thanks to the penalization by
the $\ell_1$ norm.
}
\textcolor{black}{This criterion actually boils down to the Generalized Lasso proposed by \citep{Tibshirani2011}:
\begin{equation}\label{genlasso}
 L_{\lambda}^{genlasso}(\widetilde{\bbeta})= \norm{\by-\widetilde{\bX}\widetilde{\bbeta}}_{2}^{2}+\lambda\norm{D\widetilde{\bbeta}}_{1}, \textrm{ with } \lambda>0
\end{equation}
and $D=\Sigma^{-1/2}$.}

\textcolor{black}{Since, as explained in \citep{Tibshirani2011}, some problems may occur when the rank of the design matrix is not full, we will consider in this paper the following
criterion:
\begin{equation}\label{WLasso_new}
 L_{\lambda, \eta}^{gEN}(\widetilde{\bbeta})= \norm{\by-\widetilde{\bX}\widetilde{\bbeta}}_{2}^{2}+\lambda\norm{\Sigma^{-1/2}\widetilde{\bbeta}}_{1}+\eta\norm{\widetilde{\bbeta}}_{2}^{2},\textrm{ with } \lambda, \eta>0.
\end{equation}
Since it consists in adding an $L_2$ penalty part to the Generalized Lasso as in the Elastic Net, we will call it generalized Elastic Net (gEN).
We prove in Section \ref{sec:theo} that under Assumption \ref{covX} and the Generalized Irrepresentable Condition (GIC) \eqref{GIC} given below among others,
$\widehat{\bbeta}$ is a sign-consistent estimator of $\bbeta^\star$ where $\widehat{\bbeta}$ is defined by
\begin{equation}\label{eq:genEstimator}
\widehat{\bbeta}=\bSigma^{-1/2}\widehat{\widetilde{\bbeta}},
\end{equation}
with
\begin{equation}\label{eq:genEstimatorSigma}
\widehat{\widetilde{\bbeta}}=\argmin_{\widetilde{\bbeta}}L_{\lambda, \eta}^{gEN}\left( \widetilde{\bbeta} \right),
\end{equation}
$L_{\lambda, \eta}^{gEN}\left(\widetilde{\bbeta} \right)$ being defined in Equation \eqref{WLasso_new}. 
The Generalized Irrepresentable Condition (GIC) can be stated as follows: There exist $\lambda,\eta,\alpha,\delta_4>0$ such that for all $j$,
\begin{equation}\label{GIC}
\bP\left(\left|\left((C_{21}^{n}+\frac{\eta}{n}\Sigma_{21})(C_{11}^{n}+\frac{\eta}{n}\Sigma_{11})^{-1}\left(\textrm{sign}(\bbeta_{1}^\star)+\frac{2\eta}{\lambda}\bbeta_{1}^\star\right)-\frac{2\eta}{\lambda}\Sigma_{21}\bbeta_{1}^\star\right)_j\right| \leq 1-\alpha\right)=1-o\left(e^{-n^{\delta_4}}\right).
\end{equation}
} 
\textcolor{black}{Note that GIC coincides with EIC when $\bX$ is not random and $\Sigma=\mathbb{I}_{p}$. Moreover, GIC does not
  require $C_{11}^{n}$ to be invertible. Since EIC and IC are both particular cases of GIC, if the IC or EIC holds, then there exist $\lambda$ or
$\eta$ such that the GIC holds.}

The rest of the paper is organized as follows. Section \ref{sec:theo} is devoted to the theoretical results of the paper.
\textcolor{black}{More precisely, we prove that under some mild conditions $\widehat{\bbeta}$ defined in \eqref{eq:genEstimator} is a sign-consistent
  estimator of $\bbeta^\star$. To support our theoretical results, some numerical experiments are presented in Section \ref{sec:num}.
The proofs of our theoretical results can be found in Section \ref{sec:proofs}}.



\section{Theoretical results}\label{sec:theo}

\textcolor{black}{The goal of this section is to establish the sign consistency of the Generalized Elastic Net estimator defined in (\ref{eq:genEstimator}).} 
\textcolor{black}{To prove this result, we shall use the following lemma.}
\textcolor{black}{
  \begin{lem}\label{lem1}
    Let $\by$ satisfying Model (\ref{lm}) under Assumption \ref{covX} and $\widehat{\bbeta}$ be defined in \eqref{eq:genEstimator}. Then,
\begin{equation}
\bP\left(sign\left(\widehat{\bbeta}\right)=sign(\bbeta^\star) \right) \geq \bP\left(A_n\cap B_n \right),
\end{equation}
where
\[
A_n:=\left\{ \left| \left(\mathcal{C}_{11}^{n,\bSigma}\right)^{-1}W_n(1) \right| < \sqrt{n}\left( \left|\bbeta_1^\star\right| - \frac{\lambda}{2n}\left|\left(\mathcal{C}_{11}^{n,\bSigma}\right)^{-1}sign(\bbeta_1^\star)\right| - \frac{\eta}{n}\left|\left(\mathcal{C}_{11}^{n,\bSigma}\right)^{-1}\bSigma_{11}\bbeta_1^\star\right|\right) \right\},
\]
\begin{multline*}
B_n:=\left\{ \left|\mathcal{C}_{21}^{n,\bSigma}\left(\mathcal{C}_{11}^{n,\bSigma}\right)^{-1}W_n(1)-W_n(2)\right| \leq \frac{\lambda}{2\sqrt{n}} \right.\\
 \left. -\frac{\lambda}{2\sqrt{n}}\left| \mathcal{C}_{21}^{n,\bSigma}\left(\mathcal{C}_{11}^{n,\bSigma}\right)^{-1}\left(sign(\bbeta_1^\star) +\frac{2\eta}{\lambda}  \bSigma_{11}\bbeta_1^\star \right)-\frac{2\eta}{\lambda}\bSigma_{21}\bbeta_1^\star \right| 
  \right\},
\end{multline*}
and
\begin{equation}\label{eq:C11-n,Sigma}
\mathcal{C}_{11}^{n,\bSigma}=C_{11}^n+\frac{\eta}{n}\bSigma_{11},\;\mathcal{C}_{21}^{n,\bSigma}=C_{21}^n+\frac{\eta}{n}\bSigma_{21},\;
W_n=\dfrac{1}{\sqrt{n}}\bX'\bepsilon=\begin{bmatrix}W_n(1)\\W_n(2)\end{bmatrix},
\end{equation}
with
$$
W_n(1)=\dfrac{1}{\sqrt n}\bX'_1\bepsilon \textrm{ and } W_n(2)=\dfrac{1}{\sqrt n}\bX'_{2}\bepsilon.
$$
\end{lem}
}

\textcolor{black}{The proof of Lemma \ref{lem1} is given in Section \ref{sec:proofs}.}

\textcolor{black}{The following theorem gives the conditions under which the sign consistency of the generalized Elastic Net estimator $\widehat{\bbeta}$
  defined in (\ref{eq:genEstimator}) holds.}

\textcolor{black}{
\begin{thm}\label{thm1}
  Assume that $\by$ satisfies Model (\ref{lm}) under Assumption \ref{covX}
  with $p=p_n$ is such that $p_n\exp\left(n^{-\delta}\right)$ tends to 0 as $n$ tends to infinity for all positive $\delta$.
  Assume also that there exist some positive constants $M_1$, $M_2$, $M_3$ and $\alpha$ satisfying
  \begin{equation}\label{eq:M1_M2_M3}
    M_1<\frac{\bbeta_{\min}^2}{9\sigma^2}\quad \textrm{ and }\quad \dfrac{\sqrt{2+\sqrt{2}}\sqrt{M_3}\sigma}{\alpha}<\dfrac{\bbeta_{\min}}{3M_2\sqrt{q}},
  \end{equation}
  and that there exist $\lambda> 0$ and $\eta> 0$ such that (\ref{GIC}) and
\begin{equation}\label{eq:lambda/n<}
\dfrac{\lambda}{n}<\dfrac{2\bbeta_{\min}}{3M_2\sqrt{q}},
\end{equation}
\begin{equation}\label{eq:lambda/n>}
\dfrac{\lambda}{n}\geq\dfrac{2\sqrt{2+\sqrt{2}}\sqrt{M_3}\sigma}{\alpha},
\end{equation}
\begin{equation}\label{eq:eta/n<}
\dfrac{\eta}{n}<\dfrac{1}{3M_2\lambda_{\max}\left(\bSigma_{11}\right)}\times \dfrac{\bbeta_{\min}}{\left\| \bbeta_1^\star \right\|_2},
\end{equation}
hold as $n$ tends to infinity,
where $\bbeta_{\min}=\min_{1\leq j\leq q}\left|\left(\bbeta_1^\star \right)_j\right|$.
  Suppose also that there exist some positive constants $\delta_1$, $\delta_2$, $\delta_3$ such that,
as $n\rightarrow\infty$,
\begin{equation}\label{eq:lambda_max_HA}
  \bP\left(\lambda_{\max}\left(H_AH'_A\right)\leq M_1 \right)= 1- o\left(e^{-n^{\delta_1}}\right),
\end{equation}
\begin{equation}\label{eq:lambda_max_C11-1}
  \bP\left(\lambda_{\max}\left(\left(\mathcal{C}_{11}^{n,\bSigma}\right)^{-1}\right)\leq M_2 \right)= 1- o\left(e^{-n^{\delta_2}}\right),
\end{equation}
\begin{equation}\label{eq:lambda_max_HB}
  \bP\left(\lambda_{\max}\left(H_BH'_B\right)\leq M_3 \right)= 1- o\left(e^{-n^{\delta_3}}\right),
\end{equation}
where $\lambda_{\max}(A)$ denotes the largest eigenvalue of $A$,
$$
H_A=\frac{1}{\sqrt{n}}\left(\mathcal{C}_{11}^{n,\bSigma}\right)^{-1}\bX'_1 \textrm{ and } H_B=\dfrac{1}{\sqrt{n}}\left(\mathcal{C}_{21}^{n,\bSigma}
  \left(\mathcal{C}_{11}^{n,\bSigma}\right)^{-1}\bX'_1 - \bX'_{2} \right),
$$
$\mathcal{C}_{11}^{n,\bSigma}$ and  $\mathcal{C}_{21}^{n,\bSigma}$ being defined in (\ref{eq:C11-n,Sigma})
and $\bX_1$ (resp. $\bX_2$) denoting the first $q$ (resp. the last $p-q$) columns of $\bX$.
Then,
$$
\bP\left(sign\left(\widehat{\bbeta}\right)=sign(\bbeta^\star) \right)\rightarrow 1, \textrm{ as } n\rightarrow\infty,
$$
where $\widehat{\bbeta}$ is defined in \eqref{eq:genEstimator}.
\end{thm}
Note that Conditions (\ref{eq:lambda/n<}) and (\ref{eq:lambda/n>}) are consistent thanks to (\ref{eq:M1_M2_M3}).
}

\textcolor{black}{The proof of Theorem \ref{thm1} is given in Section \ref{sec:proofs} and a discussion on the assumptions of Theorem \ref{thm1}
is provided in Section \ref{sec:num}.}


\section{Numerical experiments}\label{sec:num}

\textcolor{black}{The goal of this section is to discuss the assumptions and illustrate the results of Theorem \ref{thm1}.
  For this, we generated datasets from Model (\ref{lm}) where the matrix $\bSigma$ appearing in \ref{covX} is defined by
\begin{equation}
     \label{eq:SPAC}
     \bSigma=
       \begin{bmatrix}
         \bSigma_{11} &  \bSigma_{12} \\
         \bSigma_{12}' &  \bSigma_{22}
       \end{bmatrix}.
       \vspace{-1mm}
      \end{equation} 
      In (\ref{eq:SPAC}), $\bSigma_{11}$ is the correlation matrix of the active variables having its off-diagonal entries equal to $\alpha_1$, $\bSigma_{22}$  is the
      correlation matrix of the
      non active variables having its off-diagonal entries equal to 
      $\alpha_3$ and $\bSigma_{12}$ is the correlation matrix between the active and the non active variables with entries equal to $\alpha_2$. In the numerical
      experiments, $(\alpha_{1} , \alpha_{2} , \alpha_{3})=(0.3, 0.5, 0.7)$.
      Moreover, $\bbeta^\star$ appearing in Model (\ref{lm}) has $q$ non zero components which are equal to $b$ and $\sigma=1$.
The number of predictors $p$ is equal to 200, 400, or 600 and the sample size $n$ takes the same values for each value of $p$.}

\subsection{Discussion on the assumptions of Theorem \ref{thm1}}

\begin{figure}[!h]
  \begin{center}
    \includegraphics[width=0.8\textwidth]{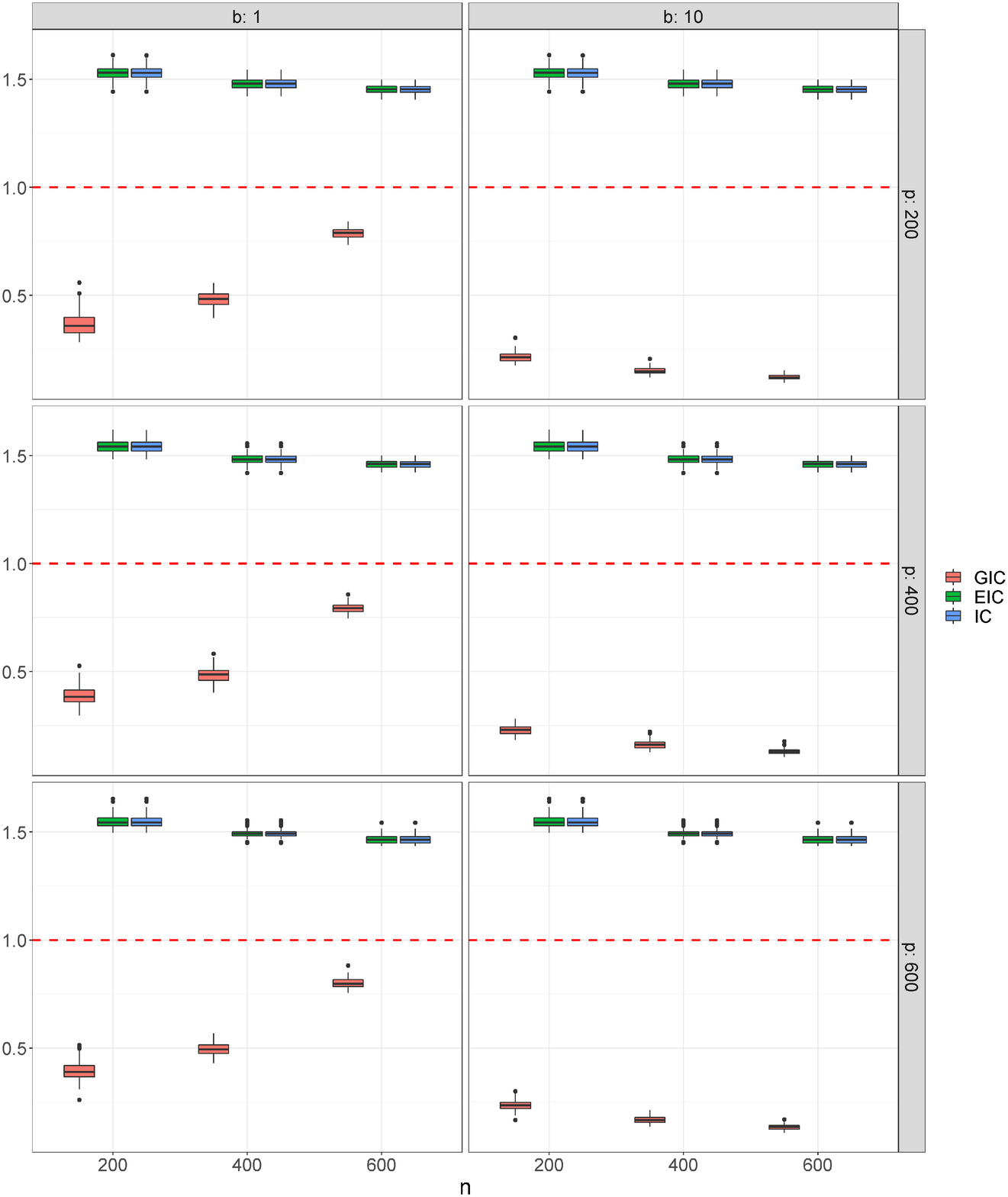}
  \end{center}
  \caption{Boxplot of values defined in (\ref{IC_value}) and obtained from 100 replications.\label{fig:IC_box}}
\end{figure}

\textcolor{black}{We first show that GIC defined in (\ref{GIC}) can be satisfied even when EIC and IC, defined in (\ref{EIC}) and (\ref{IC}) respectively,
  are not fulfilled. For this, we computed for different values of $\lambda$
and $\eta$ the following values:
\begin{align}\label{IC_value}
  &\textrm{IC}= \max_j\left(\left|\left(C_{21}^{n}(C_{11}^{n})^{-1}(\textrm{sign}(\bbeta_{1}^\star)\right)_j\right|\right)\nonumber\\
  &\textrm{EIC}= \min_{\lambda,\eta}\max_j\left(\left|\left(C_{21}^{n}(C_{11}^{n}+\frac{\eta}{n}\mathbb{I}_{q})^{-1}(\textrm{sign}(\bbeta_{1}^\star)+\frac{2\eta}{\lambda}
    \bbeta_{1}^\star)\right)_j\right|\right)\nonumber\\
 &\textrm{GIC}= \min_{\lambda,\eta}\max_j\left(\left|\left((C_{21}^{n}+\frac{\eta}{n}\Sigma_{21})(C_{11}^{n}+\frac{\eta}{n}\Sigma_{11})^{-1}\left(\textrm{sign}(\bbeta_{1}^\star)+\frac{2\eta}{\lambda}\bbeta_{1}^\star\right)-\frac{2\eta}{\lambda}\Sigma_{21}\bbeta_{1}^\star\right)_j\right|\right)
\end{align}
and Figure \ref{fig:IC_box} displays the boxplots of these criteria obtained from 100 replications.
We can see from these figures that in all the considered cases GIC is satisfied (\textit{i.e.} all values are smaller than 1) whereas EIC and IC are not.
The values of $p$ and $n$ do not seem to have a big impact on EIC and IC. However, contrary to $p$, $n$ seems to have an influence on GIC which
increases with $n$ when $b=1$ and decreases when $n$ increases when $b=10$.
}

\begin{figure}[!h]
  \begin{center}
    \includegraphics[width=0.8\textwidth]{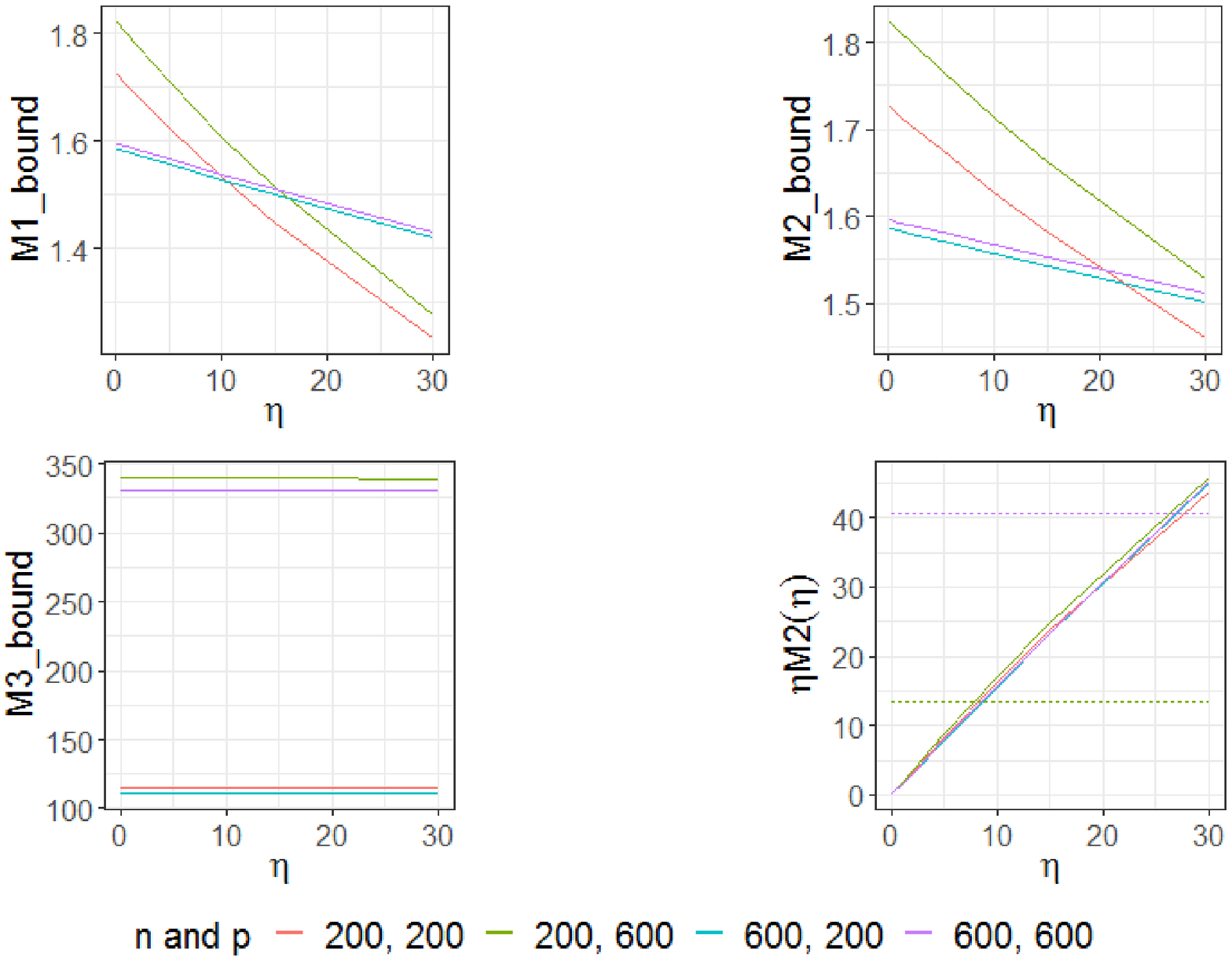}
  \end{center}
  \caption{Top left: Average of $\lambda_{\max}\left(H_AH'_A\right)$ in (\ref{eq:lambda_max_HA}) as a function of $\eta$.
    Top right: Average of $\lambda_{\max}\left(\left(\mathcal{C}_{11}^{n,\bSigma}\right)^{-1}\right)$ in (\ref{eq:lambda_max_C11-1}) as a function of $\eta$.
    Bottom left: Average of $\lambda_{\max}\left(H_BH'_B\right)$ in (\ref{eq:lambda_max_HB}) as a function of $\eta$.
    Bottom right: Average of the left (resp. right) part of (\ref{eq18_new}) in plain (resp. dashed) line. The averages are obtained from 10 replications. Here $q=5$.\label{fig:cond_q5}}
\end{figure}

\begin{figure}[!h]
  \begin{center}
    \includegraphics[width=0.8\textwidth]{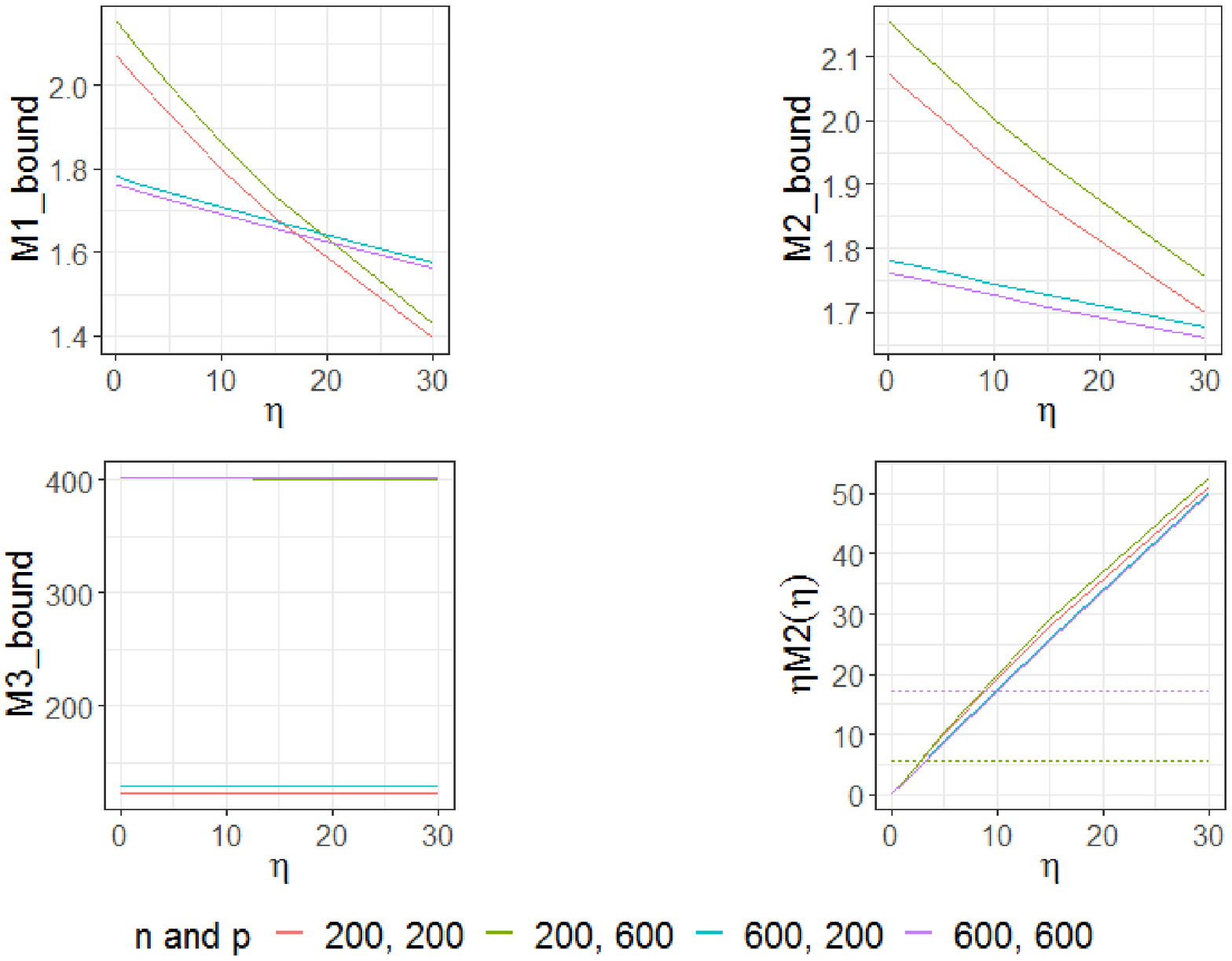}
  \end{center}
  \caption{Top left: Average of $\lambda_{\max}\left(H_AH'_A\right)$ in (\ref{eq:lambda_max_HA}) as a function of $\eta$.
    Top right: Average of $\lambda_{\max}\left(\left(\mathcal{C}_{11}^{n,\bSigma}\right)^{-1}\right)$ in (\ref{eq:lambda_max_C11-1}) as a function of $\eta$.
    Bottom left: Average of $\lambda_{\max}\left(H_BH'_B\right)$ in (\ref{eq:lambda_max_HB}) as a function of $\eta$.
    Bottom right: Average of the left (resp. right) part of (\ref{eq18_new}) in plain (resp. dashed) line. The averages are obtained from 10 replications. Here $q=10$.\label{fig:cond_q10}}
\end{figure}

\textcolor{black}{Figures \ref{fig:cond_q5} and \ref{fig:cond_q10} show  the behavior of $\lambda_{\max}\left(H_AH'_A\right)$,
  $\lambda_{\max}\left(\left(\mathcal{C}_{11}^{n,\bSigma}\right)^{-1}\right)$ and $\lambda_{\max}\left(H_BH'_B\right)$
appearing in (\ref{eq:lambda_max_HA}), (\ref{eq:lambda_max_C11-1}) and (\ref{eq:lambda_max_HB}) with respect to $\eta$ for different values of $n$, $p$
and for $q=5$ or 10. These plots thus provide lower bounds for $M_1$, $M_2$ and $M_3$ appearing in the previous equations.
Observe that (\ref{eq:eta/n<}) can be rewritten as:
\begin{equation}\label{eq18_new}
\eta M_2 < \dfrac{n}{3\lambda_{\max}\left(\bSigma_{11}\right)}\times \dfrac{\bbeta_{\min}}{\left\| \bbeta_1^\star \right\|_2}.
\end{equation}
Based on the plots at the bottom right of Figures \ref{fig:cond_q5} and \ref{fig:cond_q10}, we can see
that there exist $\eta$'s satisfying Condition \ref{eq18_new} and thus (\ref{eq:eta/n<}) and that the interval in which the adapted $\eta$'s
lie is larger when $q=5$ than when $q=10$.}

\begin{figure}[!h]
  \begin{center}
    \includegraphics[width=0.8\textwidth]{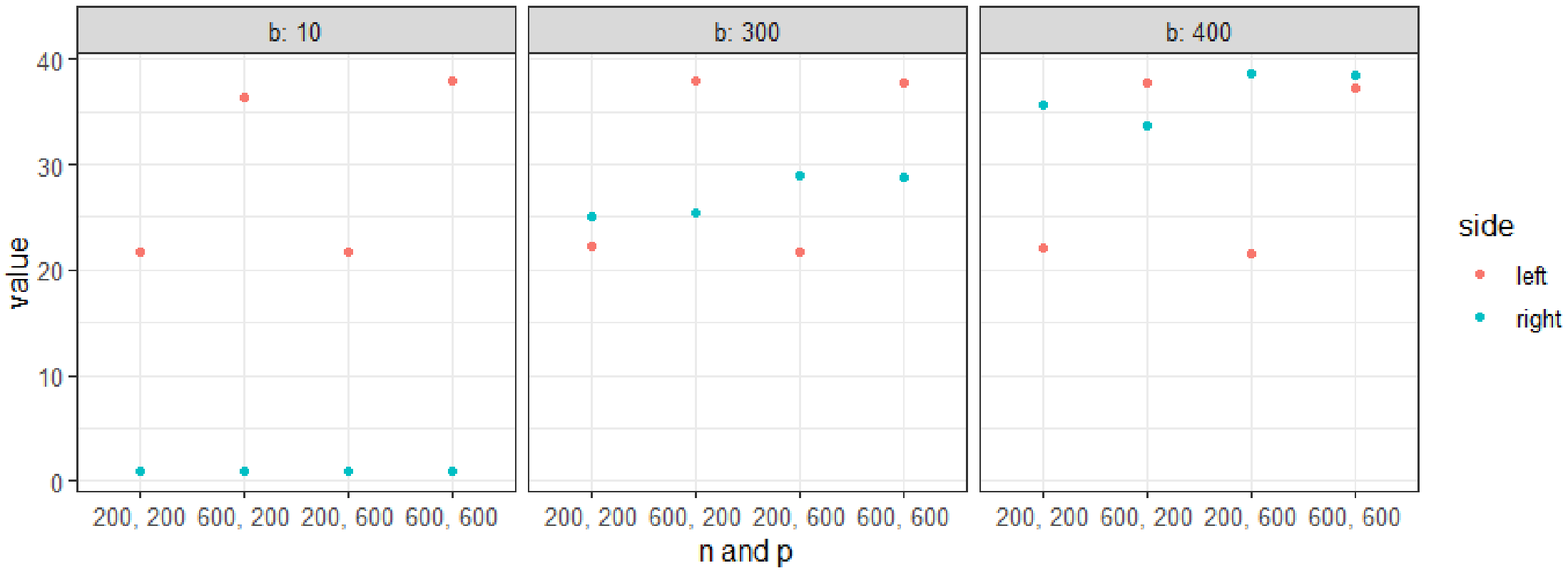}
  \end{center}
  \vspace{-25mm}
  \caption{Average of the left-hand (resp. right-hand) side of the second part of (\ref{eq:M1_M2_M3}) in red (resp. blue) for $q=5$.\label{fig:eq15_q5}}
\end{figure}

\begin{figure}[!h]
  \begin{center}
    \includegraphics[width=0.8\textwidth]{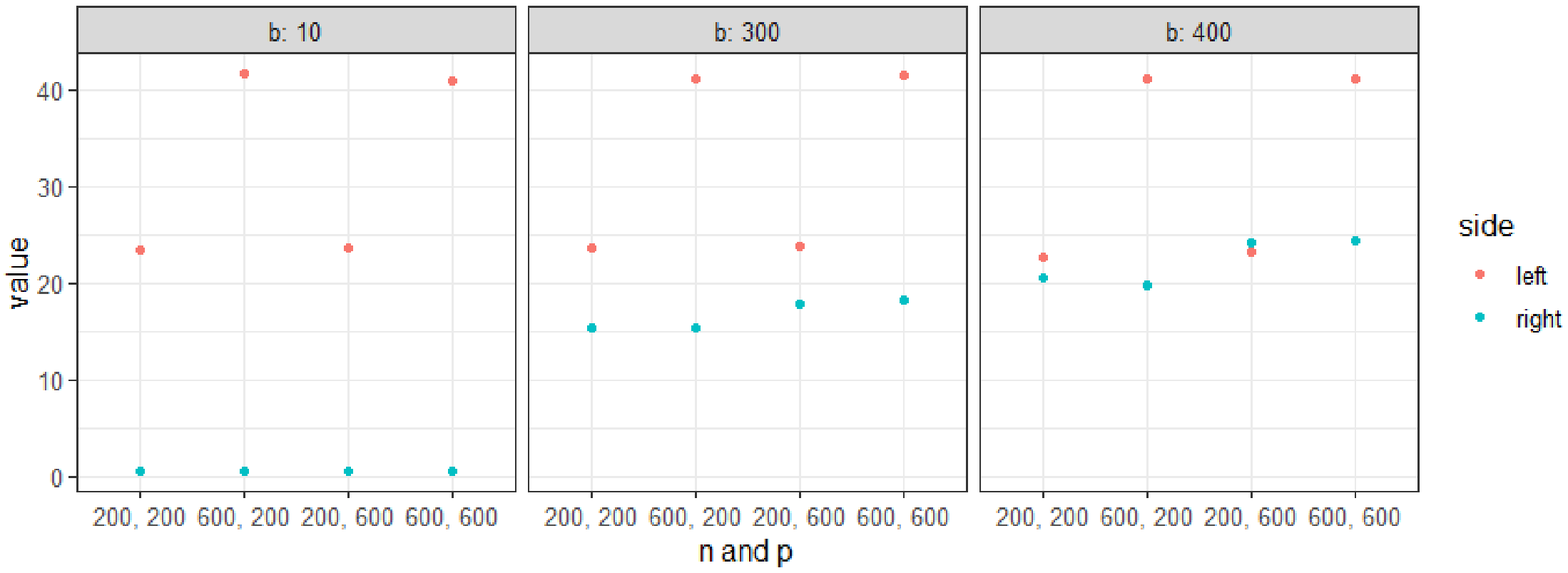}
  \end{center}
  \vspace{-25mm}
  \caption{Average of the left-hand (resp. right-hand) side of the second part of (\ref{eq:M1_M2_M3}) in red (resp. blue) for $q=10$.\label{fig:eq15_q10}}
\end{figure}

\textcolor{black}{Based on the average of $M_1$ previously obtained, the left part of (\ref{eq:M1_M2_M3})
  is always satisfied as soon as $b>\sqrt{18}$.
  Based on the average of $M_2$ and $M_3$ previously obtained, the average of left-hand side and of the right-hand side of the right part of
  Equation (\ref{eq:M1_M2_M3}) are displayed 
in Figures \ref{fig:eq15_q5} and \ref{fig:eq15_q10}. We can see from these figures that it
is only satisfied for large values of $b$. Moreover, it is more often satisfied when $q=5$ than for $q=10$.
}

\textcolor{black}{We will show in the next section that even if the cases where all the conditions of the theorem are not fulfilled our method is robust enough
  to outperform the Elastic Net defined in~(\ref{eq:def_EN}) even in these cases.}

\subsection{Comparison with other methods}

\textcolor{black}{To assess the performance of our approach (gEN) in terms of sign-consistency with respect to other methods
  and to illustrate the results of Theorem \ref{thm1},
  we computed the True Positive Rate (TPR), namely the proportion of active variables selected, and
  the False Positive Rate (FPR), namely the proportion of non active variables selected,
  of the Elastic Net and gEN estimators defined in (\ref{eq:def_EN}) and (\ref{eq:genEstimator}), respectively.}
  
\begin{figure}
\begin{center}
\includegraphics[trim={0.6cm 0 2.3cm 0 },clip,width=0.8\textwidth,angle=-90]{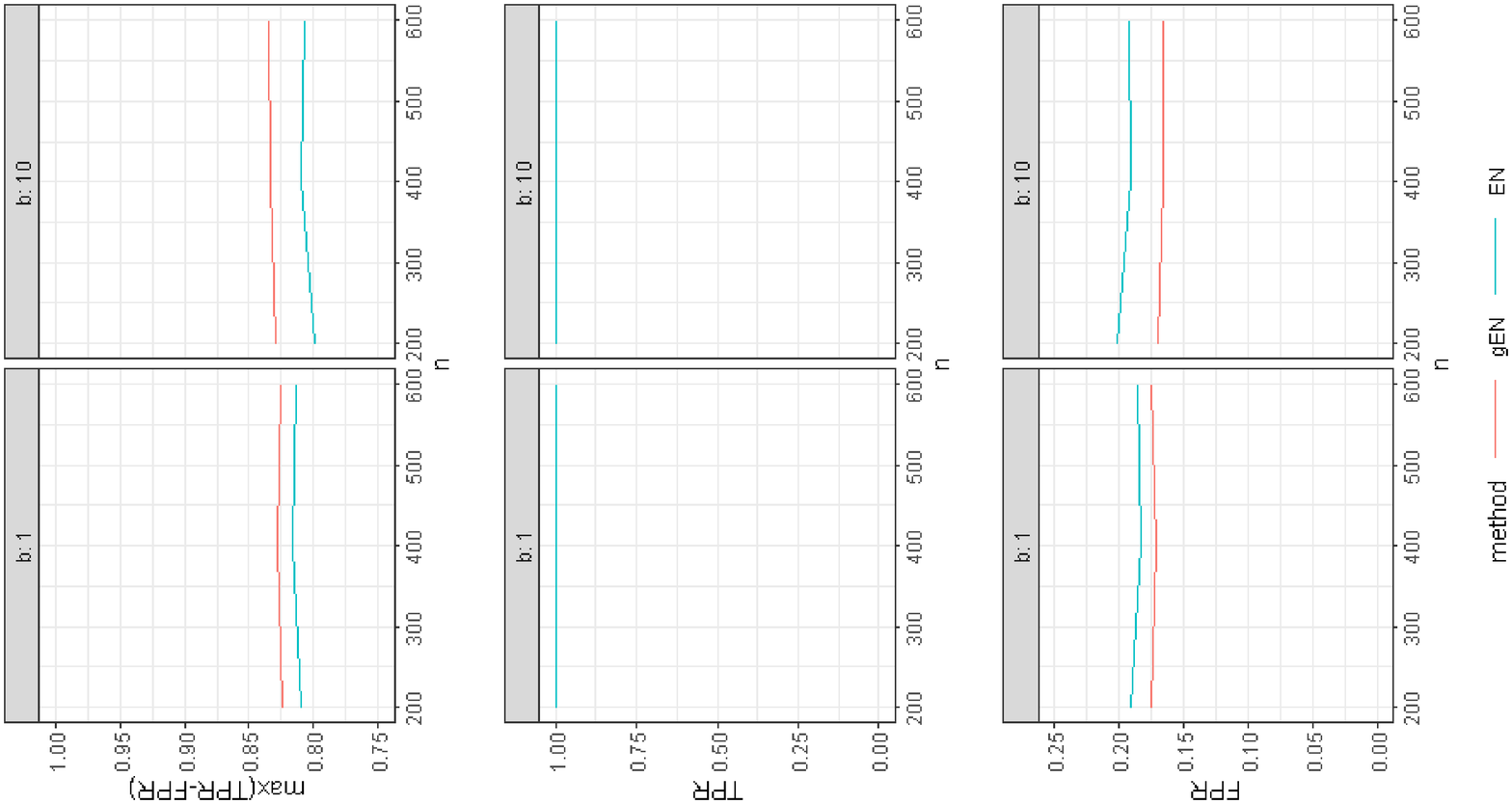}
\end{center}
\caption{Average of max(TPR-FPR) and the corresponding TPR and FPR for gEN (in red) and Elastic Net (in blue) with $p=200$.\label{fig:Res_p200}}
\end{figure}

\begin{figure}
\begin{center}
\includegraphics[trim={0.6cm 0 2.3cm 0 },clip,width=0.8\textwidth,angle=-90]{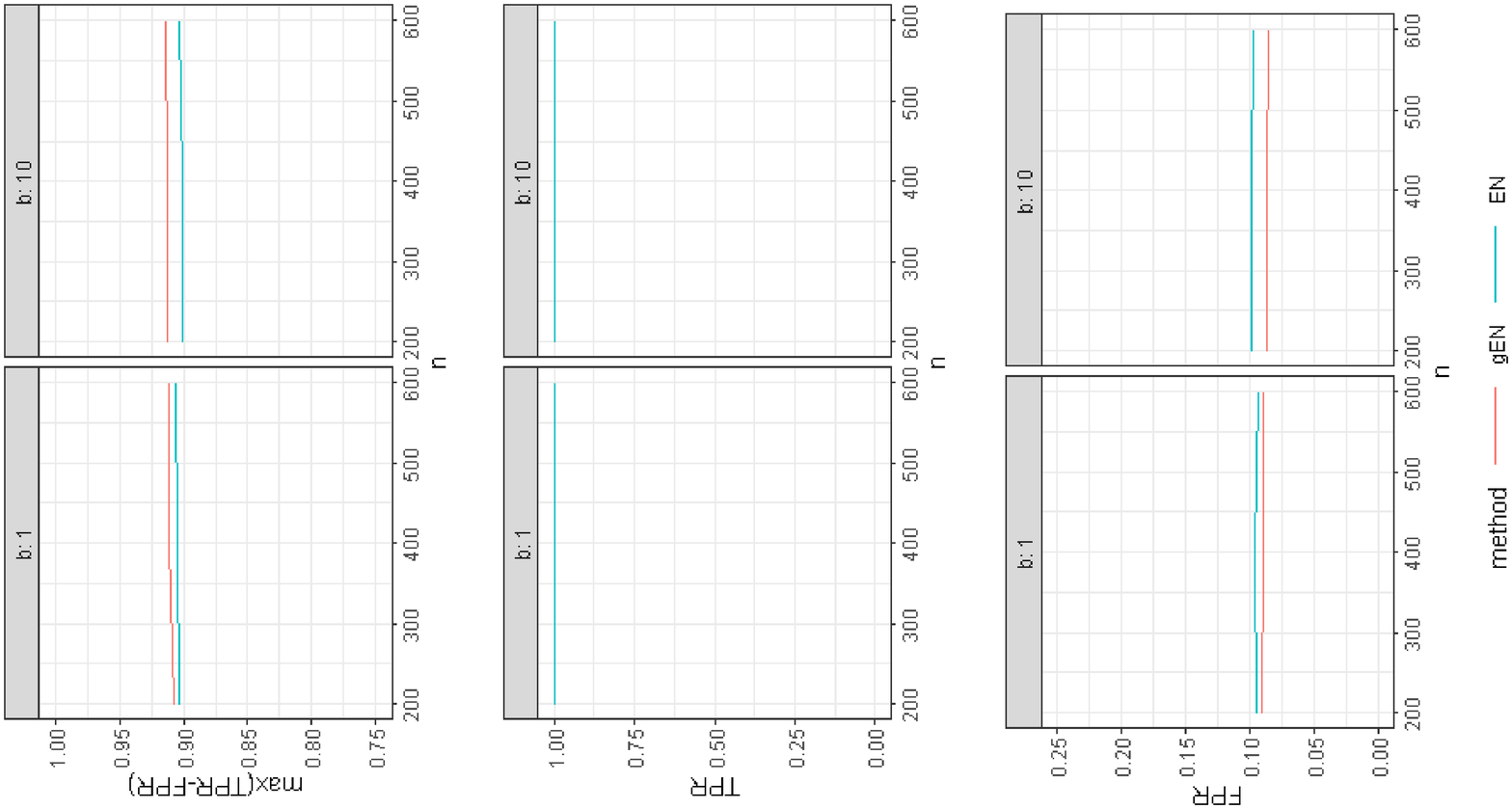}
\end{center}
\caption{Average of max(TPR-FPR) and the corresponding TPR and FPR for gEN (in red) and Elastic Net (in blue) with $p=600$.\label{fig:Res_p600}}
\end{figure}

\textcolor{black}{Figures \ref{fig:Res_p200} and \ref{fig:Res_p600} display the empirical mean
of the largest difference between the True Positive Rate and
False Positive Rate  over the replications. It is obtained by selecting for each
replication the value of $\lambda$ and $\eta$ achieving the largest difference between
the TPR and FPR and by averaging these differences. They also display
 the corresponding TPR and FPR for gEN and Elastic Net for different values of $n$ and $p$. We can see from these figures
  that the gEN and the Elastic Net estimators have a TPR equal to 1 but that the FPR of gEN is smaller than Elastic Net.
  We can see from these figures that the difference between the performance of gEN and Elastic Net is larger for high signal-to-noise ratios ($b=10$).
It has to be noticed that when TPR=1 for our approach it also means that the signs of the non null $\beta_i^\star$ are also properly retrieved.}


\begin{figure}
\begin{center}
\includegraphics[trim={0.6cm 0 3cm 0 },clip,width=0.8\textwidth,angle=-90]{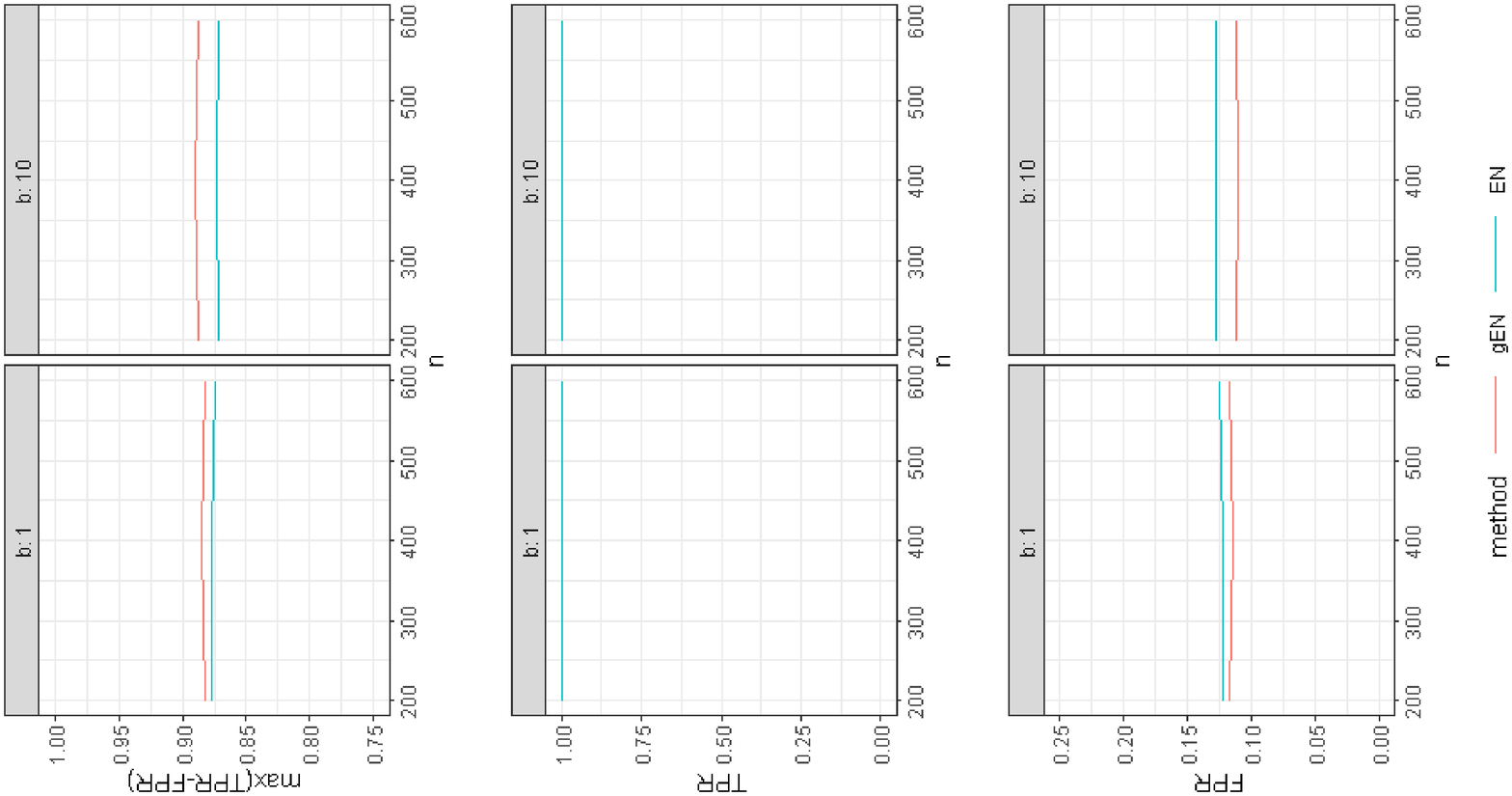}
\end{center}
\caption{Average of max(TPR-FPR) and the corresponding TPR and FPR for gEN (in red) and Elastic Net (in blue) with $p=400$.\label{fig:Res_p400}}
\end{figure}

\section{Discussion}\label{sec:discussion}
In this paper, we proposed a novel variable selection approach called gEN (generalized Elastic Net)
  in the framework of linear models where the columns
  of the design matrix are highly correlated and thus when the standard Lasso criterion usually fails.
  We proved that under mild conditions, among which the GIC, which is valid when other standard conditions like EIC or IC are not fulfilled,
  our method provides a sign-consistent estimator of $\bbeta^\star$. For a more thorough discussion regarding the application of our approach
  in practical situations, we refer the reader to \cite{zhu2020}.

\section{Proofs}\label{sec:proofs}
\subsection{Proof of Lemma \ref{lem1}}
Note that (\ref{WLasso_new}) given by:
$$
L^{gEN}(\widetilde{\bbeta})=\norm{\by-\widetilde{\bX}\widetilde{\bbeta}}_{2}^{2}+\lambda\norm{\Sigma^{-1/2}\widetilde{\bbeta}}_{1}+\eta\norm{\widetilde{\bbeta}}_{2}^{2}
$$
can be rewritten as 
\begin{equation*}
L^{gEN}(\widetilde{\bbeta})= \norm{\by^{*}-\widetilde{\bX}{*}\widetilde{\bbeta}}_{2}^{2}+\lambda\norm{\bSigma^{-1/2}\widetilde{\bbeta}}_{1}, 
\end{equation*}
where
\begin{equation*}
\by^{*} = 
\begin{pmatrix}
\by \\
0
\end{pmatrix}, \quad
\widetilde{\bX}^{*} = 
\begin{pmatrix}
\widetilde{\bX} \\
\sqrt{\eta}\mathbb{I}_{p}
\end{pmatrix}.
\end{equation*}
Then, $\widehat{\widetilde{\bbeta}}$ satisfies 
\begin{equation}\label{eq:kkt}
 \widetilde{\bX}^{*'} \left(\by^{*}-\widetilde{\bX}{*'}\widehat{\widetilde{\bbeta}}\right)=\frac{\lambda}{2}(\bSigma^{-1/2})'z,
\end{equation}
where $A'$ denotes the transpose of the matrix $A$, and 
\begin{equation*}
    \begin{cases}
      z_{j}=sign\left((\bSigma^{-1/2}\widehat{\widetilde{\bbeta}})_j\right), & \text{if}\ (\bSigma^{-1/2}\widehat{\widetilde{\bbeta}})_{j} \neq 0 \\
      z_{j}\in [-1,1], & \text{if}\ (\bSigma^{-1/2}\widehat{\widetilde{\bbeta}})_{j} = 0
    \end{cases}.
\end{equation*}
Equation~(\ref{eq:kkt}) can be rewritten as:
\[
\bX'\by-\left(\bX'\bX+\eta\bSigma\right)\widehat{\bbeta}=\frac{\lambda}{2}z
\]
which leads to
\[
\bX'\bX(\bbeta^\star-\widehat{\bbeta})+\bX'\bepsilon-\eta\bSigma\widehat{\bbeta}=\frac{\lambda}{2}z,
\]
by using that $\by=\bX\bbeta^\star+\bepsilon$.
By using the following notations: $\widehat{\bu}=\widehat{\bbeta}-\bbeta^\star$,
$$
C_n=\frac{1}{n}\bX'\bX \textrm{ and } W_n=\frac{1}{\sqrt{n}}\bX'\bepsilon,
$$
Equation~(\ref{eq:kkt}) becomes
\begin{equation}\label{eq:kkt2}
\left(C_n+\frac{\eta}{n}\bSigma\right)\sqrt{n}\widehat{\bu}+\frac{\eta}{\sqrt{n}}\bSigma\bbeta^\star-W_n=-\frac{\lambda}{2\sqrt{n}}z.
\end{equation}
With the following notations:
$$
C_n=\begin{pmatrix}C_{11}^n & C_{12}^n \\ C_{21}^n & C_{22}^n\end{pmatrix},\;
\bSigma=\begin{pmatrix}\bSigma_{11} & \bSigma_{12} \\ \bSigma_{21} & \bSigma_{22}\end{pmatrix},\;
\widehat{\bu}=\begin{pmatrix}\widehat{\bu}_1 \\ \widehat{\bu}_2\end{pmatrix},\;
W_n=\begin{pmatrix}W_n(1) \\ W_n(2)\end{pmatrix},\;
\bbeta^\star=\begin{pmatrix}\bbeta_1^\star \\ \textbf{0} \end{pmatrix},
$$
the first components of Equation~(\ref{eq:kkt2}) are:
\begin{equation}\label{eq:kkt2_1}
  \left(C_{11}^n+\frac{\eta}{n}\bSigma_{11}\right)\sqrt{n}\widehat{\bu}_1+\left(C_{12}^n+\frac{\eta}{n}\bSigma_{12}\right)\sqrt{n}\widehat{\bu}_2+
  \frac{\eta}{\sqrt{n}}\bSigma_{11}\bbeta_1^\star-W_n(1)=-\frac{\lambda}{2\sqrt{n}}sign(\bbeta_1^\star).
\end{equation}
If $\widehat{\bu}=\begin{pmatrix}\widehat{\bu}_1 \\ 0\end{pmatrix}$, it can be seen as a solution of the generalized Elastic Net criterion where, by
Equation~(\ref{eq:kkt2_1}), $\widehat{\bu}_1$ is defined by: 
\begin{equation}\label{u1_chap}
  \sqrt{n}\widehat{\bu}_1=\left(\mathcal{C}_{11}^{n,\bSigma}\right)^{-1}W_n(1)-\frac{\eta}{\sqrt{n}}\left(\mathcal{C}_{11}^{n,\bSigma}\right)^{-1}
  \bSigma_{11}{\bbeta_1}^\star-\frac{\lambda}{2\sqrt{n}}\left(\mathcal{C}_{11}^{n,\bSigma}\right)^{-1}sign({\bbeta_1}^\star),
\end{equation}
where we used (\ref{eq:C11-n,Sigma}).

\textcolor{black}{Note that the event $A_n$ can be rewritten as follows:
\begin{multline}
\sqrt{n}\left( -\left|\bbeta_1^\star\right| + \frac{\lambda}{2n}\left|\left(\mathcal{C}_{11}^{n,\bSigma}\right)^{-1}sign(\bbeta_1^\star)\right| + \frac{\eta}{n}\left|\left(\mathcal{C}_{11}^{n,\bSigma}\right)^{-1}\bSigma_{11}\bbeta_1^\star\right|\right)\\
<\left(\mathcal{C}_{11}^{n,\bSigma}\right)^{-1}W_n(1)<\\
\sqrt{n}\left( \left|\bbeta_1^\star\right| - \frac{\lambda}{2n}\left|\left(\mathcal{C}_{11}^{n,\bSigma}\right)^{-1}sign(\bbeta_1^\star)\right| - \frac{\eta}{n}\left|\left(\mathcal{C}_{11}^{n,\bSigma}\right)^{-1}\bSigma_{11}\bbeta_1^\star\right|\right)\nonumber
\end{multline}
which~implies
\begin{multline}
\sqrt{n}\left( -\left|\bbeta_1^\star\right| + \frac{\lambda}{2n}\left(\mathcal{C}_{11}^{n,\bSigma}\right)^{-1}sign(\bbeta_1^\star) + \frac{\eta}{n}\left(\mathcal{C}_{11}^{n,\bSigma}\right)^{-1}\bSigma_{11}\bbeta_1^\star\right)\\
<\left(\mathcal{C}_{11}^{n,\bSigma}\right)^{-1}W_n(1)<\\
\sqrt{n}\left( \left|\bbeta_1^\star\right| + \frac{\lambda}{2n}\left(\mathcal{C}_{11}^{n,\bSigma}\right)^{-1}sign(\bbeta_1^\star) + \frac{\eta}{n}\left(\mathcal{C}_{11}^{n,\bSigma}\right)^{-1}\bSigma_{11}\bbeta_1^\star\right),\nonumber
\end{multline}
using that $-|x|\leq x\leq |x|,\,\forall x\in\R$.
Then, by using (\ref{u1_chap}), we get that $\sqrt{n}|\widehat{\bu}_1|<\sqrt{n}|\bbeta_1^\star|$ and thus $|\widehat{\bu}_1|<|\bbeta_1^\star|$.
Notice that $|\widehat{\bu}_1|<|\bbeta_1^\star|$ implies that $\widehat{\bbeta}_1\neq 0$ and that $sign(\widehat{\bbeta}_1)=sign(\bbeta_1^\star)$.
Moreover, since $\widehat{\bu}_2=0$, we get that $sign(\widehat{\bbeta})=sign(\bbeta^\star)$.}

\textcolor{black}{The last components  of (\ref{eq:kkt2}) satisfy:
\begin{equation*}
  \left(C_{21}^n+\frac{\eta}{n}\bSigma_{21}\right)\sqrt{n}\widehat{\bu}_1+\left(C_{22}^n+\frac{\eta}{n}\bSigma_{22}\right)\sqrt{n}\widehat{\bu}_2+
  \frac{\eta}{\sqrt{n}}\bSigma_{21}\bbeta_1^\star-W_n(2)=-\frac{\lambda}{2\sqrt{n}}z_2,
\end{equation*}
where by (\ref{eq:kkt}), $|z_2|\leq 1$. Hence, 
\begin{equation*}
\left| \left(C_{21}^n+\frac{\eta}{n}\bSigma_{21} \right)\sqrt{n}\widehat{\bu}_1 + \frac{\eta}{\sqrt{n}}\bSigma_{21}\bbeta_1^\star - W_n(2) \right| \leq \frac{\lambda}{2\sqrt{n}},
\end{equation*}
which can be rewritten as follows by using (\ref{u1_chap}):
\begin{equation}\label{eq:calculus_Bn}
\left|\mathcal{C}_{21}\left(\mathcal{C}_{11}^{n,\bSigma}\right)^{-1}\left(W_n(1)-\frac{\eta}{\sqrt{n}}
  \bSigma_{11}\bbeta_1^\star-\frac{\lambda}{2\sqrt{n}}sign(\bbeta_1^\star)\right)
+ \frac{\eta}{\sqrt{n}}\bSigma_{21}\bbeta_1^\star - W_n(2) \right| \leq \frac{\lambda}{2\sqrt{n}}.
\end{equation}
When the event $B_n$ is satisfied:
\begin{multline}
  - \frac{\lambda}{2\sqrt{n}}+ \frac{\lambda}{2\sqrt{n}}\left| \mathcal{C}_{21}^{n,\bSigma}\left(\mathcal{C}_{11}^{n,\bSigma}\right)^{-1}\left(sign(\bbeta_1^\star) -\frac{2\eta}{\lambda}  \bSigma_{11}\bbeta_1^\star \right)-\frac{2\eta}{\lambda}\bSigma_{21}\bbeta_1^\star \right|\\
  \leq
\mathcal{C}_{21}^{n,\bSigma}\left(\mathcal{C}_{11}^{n,\bSigma}\right)^{-1}W_n(1)-W_n(2) \\\leq \frac{\lambda}{2\sqrt{n}} 
  -\frac{\lambda}{2\sqrt{n}}\left| \mathcal{C}_{21}^{n,\bSigma}\left(\mathcal{C}_{11}^{n,\bSigma}\right)^{-1}\left(sign(\bbeta_1^\star) +\frac{2\eta}{\lambda}  \bSigma_{11}\bbeta_1^\star \right)-\frac{2\eta}{\lambda}\bSigma_{21}\bbeta_1^\star \right|.
\end{multline}
By using that $-|x|\leq x\leq |x|$ for all $x$ in $\mathbb{R}$, we get that it implies that
$$
\left|\mathcal{C}_{21}^{n,\bSigma}\left(\mathcal{C}_{11}^{n,\bSigma}\right)^{-1}W_n(1)-W_n(2)
-\frac{\lambda}{2\sqrt{n}}\mathcal{C}_{21}^{n,\bSigma}\left(\mathcal{C}_{11}^{n,\bSigma}\right)^{-1}\left(sign(\bbeta_1^\star) +\frac{2\eta}{\lambda}  \bSigma_{11}\bbeta_1^\star \right)+\frac{\eta}{\sqrt{n}}\bSigma_{21}\bbeta_1^\star \right|\leq \frac{\lambda}{2\sqrt{n}},
$$
which corresponds to (\ref{eq:calculus_Bn}).
Thus, if $A_n$ and $B_n$ are satisfied, we get that $sign\left(\widehat{\bbeta} \right)=sign\left(\bbeta \right)$, which concludes the proof.}

\subsection{Proof of Theorem \ref{thm1}}
\textcolor{black}{
  By Lemma~\ref{lem1},
  $$
  \bP\left(sign\left(\widehat{\bbeta}\right)=sign(\bbeta^\star) \right) \geq \bP\left(A_n\cap B_n \right)\geq 1-\bP\left(A_n^c\right)-\bP\left( B_n^c \right),
  $$
  where $A_n^c$ and $B_n^c$ denote the complementary of $A_n$ and $B_n$, respectively.
  Thus, to prove the theorem it is enough to prove that $\bP\left(A_n^c\right)\to 0$ and $\bP\left(B_n^c\right)\to 0$ as $n\to\infty$.
}

\textcolor{black}{
 Recall that
  \[
A_n:=\left\{ \left| \left(\mathcal{C}_{11}^{n,\bSigma}\right)^{-1}W_n(1) \right| < \sqrt{n}\left( \left|\bbeta_1^\star\right| - \frac{\lambda}{2n}\left|\left(\mathcal{C}_{11}^{n,\bSigma}\right)^{-1}sign(\bbeta_1^\star)\right| - \frac{\eta}{n}\left|\left(\mathcal{C}_{11}^{n,\bSigma}\right)^{-1}\bSigma_{11}\bbeta_1^\star\right|\right) \right\}.
\]
Let $\boldsymbol{\zeta}$ and $\boldsymbol{\tau}$ be defined by
$$
\boldsymbol{\zeta}=\left(\mathcal{C}_{11}^{n,\bSigma}\right)^{-1}W_n(1) \textrm{ and }
\boldsymbol{\tau}= \sqrt{n}\left( \left|\bbeta_1^\star\right| - \frac{\lambda}{2n}\left|\left(\mathcal{C}_{11}^{n,\bSigma}\right)^{-1}sign(\bbeta_1^\star)\right| - \frac{\eta}{n}\left|\left(\mathcal{C}_{11}^{n,\bSigma}\right)^{-1}\bSigma_{11}\bbeta_1^\star\right|\right) .
$$
Then,
$$
\bP(A_n)=\bP\left(\forall j, |\zeta_j|<\tau_j\right).
$$
Thus,
$$\bP(A_n^c)=\bP\left(\exists j, |\zeta_j|\geq\tau_j\right) \leq \sum_{j=1}^q \bP\left(|\zeta_j|\geq\tau_j\right).$$
Note that
\begin{align}\label{eq:zeta_j}
\bP(|\zeta_j|\geq\tau_j)&=\bP\left(|\zeta_j|\geq \sqrt{n}\left( \left|(\bbeta_1^\star)_j\right| - \frac{\lambda}{2n}\left|\left(\left(\mathcal{C}_{11}^{n,\bSigma}\right)^{-1}sign(\bbeta_1^\star)\right)_j\right| - \frac{\eta}{n}\left|\left(\left(\mathcal{C}_{11}^{n,\bSigma}\right)^{-1}\bSigma_{11}\bbeta_1^\star\right)_j\right|\right)  \right)\nonumber\\
&= \bP\left(|\zeta_j| + \frac{\lambda}{2\sqrt{n}}\left|\left(\left(\mathcal{C}_{11}^{n,\bSigma}\right)^{-1}sign(\bbeta_1^\star)\right)_j\right| +\frac{\eta}{\sqrt{n}}\left|\left(\left(\mathcal{C}_{11}^{n,\bSigma}\right)^{-1}\bSigma_{11}\bbeta_1^\star\right)_j\right|\geq \sqrt{n} \left|(\bbeta_1^\star)_j\right|  \right)\nonumber\\
&\leq  \bP\left(|\zeta_j|\geq\sqrt{n}\frac{ \left|(\bbeta_1^\star)_j\right|}{3}\right) + \bP\left(\frac{\lambda}{2\sqrt{n}}\left|\left(\left(\mathcal{C}_{11}^{n,\bSigma}\right)^{-1}sign(\bbeta_1^\star)\right)_j\right| \geq\sqrt{n}\frac{ \left|(\bbeta_1^\star)_j\right|}{3}\right)\nonumber\\ 
&+ \bP\left(\frac{\eta}{\sqrt{n}}\left|\left(\left(\mathcal{C}_{11}^{n,\bSigma}\right)^{-1}\bSigma_{11}\bbeta_1^\star\right)_j\right|\geq\sqrt{n}\frac{ \left|(\bbeta_1^\star)_j\right|}{3}\right).
\end{align}
Observe that
$$
\boldsymbol{\zeta}=\left(\mathcal{C}_{11}^{n,\bSigma}\right)^{-1}W_n(1)=\frac{1}{\sqrt{n}}\left(C_{11}^n+\frac{\eta}{n}\bSigma_{11}\right)^{-1}\bX'_1\bepsilon = H_A\bepsilon,
$$
where
$$
H_A=\frac{1}{\sqrt{n}}\left(C_{11}^n+\frac{\eta}{n}\bSigma_{11}\right)^{-1}\bX'_1,
$$
$\bX_1$ denoting the columns of the design matrix $\bX$ associated to the $q$ active covariates.
Thus, for all $j$ in $\{1,\dots,q\}$,
$$
\zeta_j=\sum_{k=1}^n (H_A)_{jk}\bepsilon_k.
$$
By using the Cauchy-Schwarz inequality,
\begin{eqnarray*}
|\zeta_j|=\left|\sum_{k=1}^n (H_A)_{jk}\bepsilon_k\right| &\leq & \left(\sum_{k=1}^n (H_A)_{jk}^2 \right)^{1/2} \left(\sum_{k=1}^n\bepsilon_k^2 \right)^{1/2}\\
&= & \sqrt{\left(H_AH'_A\right)_{jj}}\times\|\bepsilon\|_2\\
&\leq & \sqrt{\lambda_{\max}\left(H_AH'_A\right)}\times\|\bepsilon\|_2.
\end{eqnarray*}
Hence, the first term in the r.h.s. of (\ref{eq:zeta_j}) satisfies the following inequalities: 
\begin{eqnarray}
  \bP\left(|\zeta_j|\geq\sqrt{n}\frac{ \left|(\bbeta_1^\star)_j\right|}{3}\right)
  & \leq & \bP\left(\sqrt{\lambda_{\max}\left(H_AH'_A\right)}\times\|\bepsilon\|_2\geq\sqrt{n}\frac{ \left|(\bbeta_1^\star)_j\right|}{3}\right)\nonumber\\
& \leq & \bP\left(\lambda_{\max}\left(H_AH'_A\right)\times\|\bepsilon\|_2^2 \geq n\frac{ (\bbeta_1^\star)_j^2}{9}\right).\label{eq:Ac1}
\end{eqnarray}
Since by (\ref{eq:lambda_max_HA}), there exist $M_1>0$ and $\delta_1>0$ such that
$$
\bP\left(\lambda_{\max}\left(H_AH'_A\right)\leq M_1 \right)= 1- o\left(e^{-n^{\delta_1}}\right), \textrm{ as } n\rightarrow\infty,
$$
we have:
\begin{eqnarray*}
  \bP\left(\lambda_{\max}\left(H_AH'_A\right)\times\|\bepsilon\|_2^2 \geq n\frac{ (\bbeta_1^\star)_j^2}{9}\right)
& \leq & \bP\left(\|\bepsilon\|_2^2\geq n\frac{ (\bbeta_1^\star)_j^2}{9M_1}\right) + \bP\left( \lambda_{\max}\left(H_AH'_A\right)> M_1\right) \\
& \leq & \bP\left(\dfrac{\|\bepsilon\|_2^2}{\sigma^2}\geq \dfrac{n \bbeta_{\min}^2}{9M_1\sigma^2}\right) +o\left(e^{-n^{\delta_1}}\right).
\end{eqnarray*}
Using that $\dfrac{\|\bepsilon\|_2^2}{\sigma^2}\sim\chi^2(n)$, we get, by Lemma~1 of \cite{laurent2000adaptive}, that
\begin{equation}\label{eq:Ac2}
  \bP\left(\lambda_{\max}\left(H_AH'_A\right)\times\|\bepsilon\|_2^2 \geq n\frac{ (\bbeta_1^\star)_j^2}{9}\right)
  \leq \exp\left(-\dfrac{t}{2}+\dfrac{1}{2}\sqrt{n\left(2t-n\right)} \right) +o\left(e^{-n^{\delta_1}}\right),
\end{equation}
since $t=\frac{n \bbeta_{\min}^2}{9M_1\sigma^2}>\dfrac{n}{2}$ using that $\frac{2\bbeta_{\min}^2}{9\sigma^2}>M_1$ by (\ref{eq:M1_M2_M3}).\\
By putting together Equations~(\ref{eq:Ac1}) and (\ref{eq:Ac2}) we get
\begin{equation}\label{eq:A1}
\bP\left(|\zeta_j|>\sqrt{n}\frac{ \left|(\bbeta_1^\star)_j\right|}{3}\right) \leq \exp\left(-\dfrac{t}{2}+\dfrac{1}{2}\sqrt{n\left(2t-n\right)} \right) + o\left(e^{-n^{\delta_1}}\right),
\end{equation}
with $t=\frac{n \bbeta_{\min}^2}{9M_1\sigma^2}>\dfrac{n}{2}$.\\
Let us now derive an upper bound for the second term in the r.h.s. of (\ref{eq:zeta_j}):
$$
\bP\left(\frac{\lambda}{2\sqrt{n}}\left|\left(\left(\mathcal{C}_{11}^{n,\bSigma}\right)^{-1}sign(\bbeta_1^\star)\right)_j\right| \geq\sqrt{n}\frac{ \left|(\bbeta_1^\star)_j\right|}{3}\right).
$$
By using the Cauchy-Schwarz inequality, we get that:
\begin{eqnarray*}
&&\left|\left(\left(\mathcal{C}_{11}^{n,\bSigma}\right)^{-1}sign(\bbeta_1^\star)\right)_j\right|  =  \left| \sum_{k=1}^q \left(\left(\mathcal{C}_{11}^{n,\bSigma}\right)^{-1}\right)_{jk} \left(sign(\bbeta_1^\star)\right)_k\right|\\
&\leq& \sqrt{\sum_{k=1}^q \left(\left(\mathcal{C}_{11}^{n,\bSigma}\right)^{-1}\right)_{jk}^2} \times \|sign(\bbeta_1^\star)\|_2
  \leq \sqrt{\left(\mathcal{C}_{11}^{n,\bSigma}\right)_{jj}^{-2}} \times \sqrt{q}\\
         &\leq&  \lambda_{\max}\left(\left(\mathcal{C}_{11}^{n,\bSigma}\right)^{-1} \right) \times \sqrt{q}.
\end{eqnarray*}
Then,
\begin{eqnarray}
  &&\bP\left(\frac{\lambda}{2\sqrt{n}}\left|\left(\left(\mathcal{C}_{11}^{n,\bSigma}\right)^{-1}sign(\bbeta_1^\star)\right)_j\right| \geq\sqrt{n}\frac{ \left|(\bbeta_1^\star)_j\right|}{3}\right) \nonumber\\
  &\leq&  \bP\left(\dfrac{\lambda}{2}\sqrt{\dfrac{q}{n}}\lambda_{\max}\left(\left(\mathcal{C}_{11}^{n,\bSigma}\right)^{-1} \right) \geq \sqrt{n}\frac{ \left|(\bbeta_1^\star)_j\right|}{3} \right)\nonumber\\
  &\leq & \bP\left( \lambda_{\max}\left(\left(\mathcal{C}_{11}^{n,\bSigma}\right)^{-1} \right) \geq \dfrac{2n}{3\lambda\sqrt{q}}\left|(\bbeta_1^\star)_j\right| \right)
  \nonumber\\
&\leq &\bP\left( \lambda_{\max}\left(\left(\mathcal{C}_{11}^{n,\bSigma}\right)^{-1} \right) \geq \dfrac{2n}{3\lambda\sqrt{q}}\bbeta_{\min} \right)
        =o\left(e^{-n^{\delta_2}}\right), \textrm{ as } n\to\infty,\label{eq:A2} 
\end{eqnarray}
since $\dfrac{2n}{3\lambda\sqrt{q}}\bbeta_{\min}>M_2$ by (\ref{eq:lambda/n<}).
Let us now derive an upper bound for the third term in the r.h.s. of (\ref{eq:zeta_j}):
$$\bP\left(\frac{\eta}{\sqrt{n}}\left|\left(\left(\mathcal{C}_{11}^{n,\bSigma}\right)^{-1}\bSigma_{11}\bbeta_1^\star\right)_j\right|>\sqrt{n}\frac{ \left|(\bbeta_1^\star)_j\right|}{3}\right).$$
We have
\begin{eqnarray*}
&&\left|\left(\left(\mathcal{C}_{11}^{n,\bSigma}\right)^{-1}\bSigma_{11}\bbeta_1^\star\right)_j\right|=  \left| \sum_{k=1}^q\left( \left(\mathcal{C}_{11}^{n,\bSigma}\right)^{-1}\bSigma_{11} \right)_{jk} \left( \bbeta_1^\star \right)_k \right|
\leq  \sqrt{\sum_{k=1}^q\left( \left(\mathcal{C}_{11}^{n,\bSigma}\right)^{-1}\bSigma_{11} \right)_{jk}^2}\times\left\| \bbeta_1^\star \right\|_2\\
&\leq & \sqrt{\lambda_{\max}\left(\left(\mathcal{C}_{11}^{n,\bSigma}\right)^{-1}\bSigma_{11}^2 \left(\mathcal{C}_{11}^{n,\bSigma}\right)^{-1}\right)}\times\left\| \bbeta_1^\star \right\|_2
\leq  \lambda_{\max}\left(\left(\mathcal{C}_{11}^{n,\bSigma}\right)^{-1} \right) \lambda_{\max}\left(\bSigma_{11} \right) \times\left\| \bbeta_1^\star \right\|_2.
\end{eqnarray*}
Thus,
\begin{eqnarray}
  &&\bP\left(\frac{\eta}{\sqrt{n}}\left|\left(\left(\mathcal{C}_{11}^{n,\bSigma}\right)^{-1}\bSigma_{11}\bbeta_1^\star\right)_j\right|\geq\sqrt{n}\frac{ \left|(\bbeta_1^\star)_j\right|}{3}\right)\nonumber\\
     &\leq&  \bP\left(\frac{\eta}{\sqrt{n}}\lambda_{\max}\left(\left(\mathcal{C}_{11}^{n,\bSigma}\right)^{-1} \right) \lambda_{\max}\left(\bSigma_{11} \right) \left\| \bbeta_1^\star \right\|_2 \geq  \sqrt{n}\frac{ \left|(\bbeta_1^\star)_j\right|}{3} \right)\nonumber\\
&\leq & \bP\left( \lambda_{\max}\left(\left(\mathcal{C}_{11}^{n,\bSigma}\right)^{-1} \right) \geq\dfrac{n\bbeta_{\min}}{3\eta\left\| \bbeta_1^\star \right\|_2\lambda_{\max}\left(\bSigma_{11} \right)} \right)=o\left(e^{-n^{\delta_2}}\right), \textrm{ as } n\to\infty,\label{eq:A3}
\end{eqnarray}
since $\dfrac{n\bbeta_{\min}}{3\eta\left\| \bbeta_1^\star \right\|_2\lambda_{\max}\left(\bSigma_{11} \right)}>M_2$ by (\ref{eq:eta/n<}).\\
By putting together Equations (\ref{eq:A1}), (\ref{eq:A2}) and (\ref{eq:A3}), we get:
\begin{eqnarray}
\bP\left(A_n^c\right) 
                        & \leq & q\exp\left[ -\dfrac{n}{2}\left(\kappa-\sqrt{2\kappa-1} \right) \right]+
                                 q\times o\left(e^{-n^{\delta_1}}\right) + 2q\times o\left(e^{-n^{\delta_2}}\right), \label{eq:A4}
\end{eqnarray}
with $\kappa=\frac{ \bbeta_{\min}^2}{9M_1\sigma^2}$.
Note that $\kappa-\sqrt{2\kappa-1}>0$ since $\kappa=\frac{ \bbeta_{\min}^2}{9M_1\sigma^2}> 1$ by (\ref{eq:M1_M2_M3}).
Equation~(\ref{eq:A4}) then implies that
$$
\bP\left(A_n^c\right)\rightarrow 0 \textrm{ as } n\rightarrow\infty.
$$
Let us now prove that
$
\bP\left(B_n^c\right)\rightarrow 0 \textrm{ as } n\rightarrow\infty.
$\\
Recall that
\begin{multline*}
B_n:=\left\{ \left|\mathcal{C}_{21}^{n,\bSigma}\left(\mathcal{C}_{11}^{n,\bSigma}\right)^{-1}W_n(1)-W_n(2)\right| \leq \frac{\lambda}{2\sqrt{n}} \right.\\
 \left. -\frac{\lambda}{2\sqrt{n}}\left| \mathcal{C}_{21}^{n,\bSigma}\left(\mathcal{C}_{11}^{n,\bSigma}\right)^{-1}\left(sign(\bbeta_1^\star) +\frac{2\eta}{\lambda}  \bSigma_{11}\bbeta_1^\star \right)-\frac{2\eta}{\lambda}\bSigma_{21}\bbeta_1^\star \right| 
  \right\}.
\end{multline*}
Let
$$
\boldsymbol{\psi}=\mathcal{C}_{21}^{n,\bSigma}\left(\mathcal{C}_{11}^{n,\bSigma}\right)^{-1}W_n(1)-W_n(2)
=\frac{1}{\sqrt{n}}\left(\mathcal{C}_{21}^{n,\bSigma}\left(\mathcal{C}_{11}^{n,\bSigma}\right)^{-1}\bX'_1-\bX'_2\right)\bepsilon=:H_B\bepsilon
$$
and
$$
\boldsymbol{\mu}=\frac{\lambda}{2\sqrt{n}}-
\frac{\lambda}{2\sqrt{n}}\left| \mathcal{C}_{21}^{n,\bSigma}\left(\mathcal{C}_{11}^{n,\bSigma}\right)^{-1}\left(sign(\bbeta_1^\star) +\frac{2\eta}{\lambda}  \bSigma_{11}\bbeta_1^\star \right)-\frac{2\eta}{\lambda}\bSigma_{21}\bbeta_1^\star \right|. 
$$
Then,
$$
\bP(B_n^c)=\bP(\exists j,|\psi_j|>\mu_j)\leq\sum_{j=q+1}^{p} \bP(|\psi_j|>\mu_j).
$$
By using the Cauchy-Schwarz inequality, we get that: 
 \begin{equation}\label{eq:Bc1}
 |\psi_j| = \left|\sum_{k=1}^n(H_B)_{jk}\bepsilon_k\right|
   \leq  \left(\sum_{k=1}^n(H_B)_{jk}^2 \right)^{1/2}\times \|\bepsilon\|_2
=   \sqrt{\left(H_BH'_B \right)_{jj}} \times \|\bepsilon\|_2 
\leq  \sqrt{\lambda_{\max}(H_BH'_B)} \times \|\bepsilon\|_2,
 \end{equation}
 where
 $$
 H_BH'_B=\mathcal{C}_{21}^{n,\bSigma}\left(\mathcal{C}_{11}^{n,\bSigma}\right)^{-1}C_{11}^n\left(\mathcal{C}_{11}^{n,\bSigma}\right)^{-1}\mathcal{C}_{12}^{n,\bSigma} - \mathcal{C}_{21}^{n,\bSigma}\left(\mathcal{C}_{11}^{n,\bSigma}\right)^{-1}C_{12}^n -C_{21}^n\left(\mathcal{C}_{11}^{n,\bSigma}\right)^{-1}\mathcal{C}_{12}^{n,\bSigma} +C_{22}^n.
 $$
 By (\ref{eq:lambda_max_HB}), there exist $M_3>0$ and $\delta_3>0$ such that
 $$
 \bP\left(\lambda_{\max}\left(H_BH'_B\right)\leq M_3 \right)=1-o\left(e^{-n^{\delta_3}}\right),\textrm{ as } n\to\infty.
 $$
 By the GIC condition (\ref{GIC}), there exist $\alpha>0$  and $\delta_4>0$ such that for all $j$,
 $$
 \bP\left( \left| \mathcal{C}_{21}^{n,\bSigma}\left(\mathcal{C}_{11}^{n,\bSigma}\right)^{-1}\left(sign(\bbeta_1^\star) +\frac{2\eta}{\lambda}  \bSigma_{11}\bbeta_1^\star \right)-\frac{2\eta}{\lambda}\bSigma_{21}\bbeta_1^\star \right|\leq 1-\alpha \right)=1-o\left(e^{-n^{\delta_4}}\right).
 $$
 Thus, we get that:
\begin{eqnarray}\label{eq:Bc2}
\bP(B_n^c) & \leq & \sum_{j=q+1}^{p}\bP\left(|\psi_j|>\mu_j\right) \nonumber\\
& \leq & \sum_{j=q+1}^{p} \bP\left(|\psi_j|>\dfrac{\lambda\alpha}{2\sqrt{n}}\right) + (p-q)o\left(e^{-n^{\delta_4}}\right) \nonumber\\
& \leq & \sum_{j=q+1}^{p} \bP\left(\sqrt{\lambda_{\max}(H_BH'_B)} \times \|\bepsilon\|_2>\dfrac{\lambda\alpha}{2\sqrt{n}}\right)+ (p-q)o\left(e^{-n^{\delta_4}}\right),\,\text{using Equation~(\ref{eq:Bc1})}\nonumber\\ 
& \leq & \sum_{j=q+1}^{p}\bP\left(\lambda_{\max}(H_BH'_B) \times \|\bepsilon\|_2^2>\dfrac{\lambda^2\alpha^2}{4n}\right)+(p-q)o\left(e^{-n^{\delta_4}}\right)\nonumber\\ 
& \leq & \sum_{j=q+1}^{p}\bP\left( \dfrac{\|\bepsilon\|_2^2}{\sigma^2}>\dfrac{\lambda^2\alpha^2}{4nM_3\sigma^2}\right)+(p-q)o\left(e^{-n^{\delta_3}}\right)+(p-q)o\left(e^{-n^{\delta_4}}\right) \nonumber\\ 
& \leq & (p-q)\exp\left(-\dfrac{s}{2}+\dfrac{1}{2}\sqrt{n(2s-n)} \right) +(p-q)o\left(e^{-n^{\delta_3}}\right)+(p-q)o\left(e^{-n^{\delta_4}}\right)\nonumber\\
& \leq & (p-q)\exp\left(-\dfrac{n}{2}\left(\dfrac{s}{n}-\sqrt{2\dfrac{s}{n}-1} \right)\right) +(p-q)o\left(e^{-n^{\delta_3}}\right)
         +(p-q)o\left(e^{-n^{\delta_4}}\right)\nonumber\\
&\leq & (p-q)\exp\left(-\dfrac{n}{2}\right)+(p-q)o\left(e^{-n^{\delta_3}}\right)
         +(p-q)o\left(e^{-n^{\delta_4}}\right)
\end{eqnarray}
with $\frac{s}{n}=\dfrac{\lambda^2\alpha^2}{4n^2 M_3\sigma^2}$
since $\dfrac{\lambda^2\alpha^2}{4n^2 M_3\sigma^2}\geq 2+\sqrt{2}$ by (\ref{eq:lambda/n>}).
}
\textcolor{black}{Finally, Equation~(\ref{eq:Bc2}) implies that
$$
\bP(B_n^c)\rightarrow 0, \textrm{ as } n\rightarrow\infty,
$$
which concludes the proof.
}


\bibliographystyle{chicago}

\end{document}